\documentclass[3p,times]{elsarticle}

\usepackage[english]{babel}
\usepackage{amsthm,amssymb,amsmath,amsfonts}
\usepackage{graphicx}
\usepackage{hhline,array}
\usepackage{color}
\newtheorem{lem}{\bf Lemma }[section]
\newtheorem{theorem}{\bf Theorem}[section]
\newtheorem{definition}{\bf Definition}[section]

\usepackage{ecrc}

\volume{00}

\firstpage{1}

\journalname{APNUM}

\runauth{V. S. Sizikov, D. N. Sidorov}

\jid{APNUM}

\jnltitlelogo{APNUM}

\usepackage{amssymb}

\usepackage[figuresright]{rotating}

\begin{document}

\begin{frontmatter}



\dochead{}

\title{Generalized quadrature for solving singular integral equations of Abel type 
in application to infrared tomography}

\author[label1]{V. S. Sizikov}
\author[label2]{D. N. Sidorov}

\address[label1]{ITMO University, Kronverksky pr. 49, 197101 Saint-Petersburg, Russia}
\address[label2]{Energy Systems Institute of RAS, Lermontov Str. 130, 664033 Irkutsk, Russia}

\begin{abstract}
{We propose the generalized quadrature methods for numerical solution of singular  integral equation of Abel type. We overcome the singularity using the  analytic computation of the singular integral. The problem of solution of singular integral equation is reduced to nonsingular
system of linear algebraic equations without shift meshes techniques employment.
We also propose generalized quadrature method for solution of Abel equation using 
the singular integral. Relaxed errors bounds are derived. In order to improve the accuracy  we
use 
Tikhonov regularization method. We demonstrate the efficiency of proposed techniques 
on infrared tomography problem. Numerical experiments show that it make sense to apply 
regularization in case of 
highly noisy (about 10\%) sources only. That is due to the fact that singular integral equations
enjoy selfregularization property.  }

\end{abstract}


\end{frontmatter}



\section{Introduction}

Numerical methods for solving a variety of singular integral equations (SIE) are offered in many publications, here readers may refer to \cite{lit2}--\cite{lit5},\cite{lit13,lit14,lit18,lit21,lit22,lit24,lit27,lit31,lit34,lit37,lit40,lit41}  and others. A one-dimensional SIE of the 1st and  and 2nd kind with Cauchy kernels Hilbert kernels,  logarithmic et al., as well as two-dimensional, nonlinear SIE have been addressed. 
In present article we concentrate on  Abel singular integral equation
\cite{lit1,lit2,lit6,lit7,lit8,lit9,lit14,lit22,lit24,lit26,lit35,lit36,lit37,lit39,lit40}
\begin{equation}
2\int\limits_x^R \frac{r}{\sqrt{r^2 - x^2}} k(r) \, dr = q(x),\,\,\, 0\leq x\leq R,
\label{eq1}
\end{equation}
where $k(r)$ is desired function, $q(x)$ is the source function.
Equations \eqref{eq1} are widely used in practical models including plasma diagnostics,
thermal tomography, X-ray CT, spectroscopy, galaxy clusters astrophysics, etc.
In all these problems the object of interest enjoy the axial (or spherical) symmetry.
Abel equation is also can be written as 
\begin{equation}
\int\limits_0^x \frac{k(r)}{\sqrt{x-r}}\, dr = q(x), \, 0\leq x \leq R.
\label{eq2}
\end{equation}
It has been studied in this form is \cite{lit2,lit5,lit13,lit14,lit18,lit21,lit22,lit27,
lit29,lit31,lit34,lit39,lit41}.
Equation \eqref{eq2} describes various problems in mechanics (such as tautohron problem),
scattering  and other problems.
Of course, one may transfer SIE \eqref{eq1} into SIE \eqref{eq2}
and vice versa, but it makes it more complicated to  analyse their physical meaning.

Let us below outline  the main algorithms for numerical solution of SIE
and singular integrals computation. For more details readers may refer to
\cite{lit36, lit37}.

%

\begin{enumerate}
\item
{\it Algorithms based on  relevant mesh shift.}
In papers \cite{lit3, lit4} the discrete meshes of knots with respect to  variables $r$ and $x$
are introduced, i.e. $r_j=jh, \, x_i = r_i+\Delta,$ $j,i=0,1,\dots, n,$
$r_n=R,$ where step $h=R/n,$ $\Delta$ is mesh shift which is $h/2$ \cite{lit3}
or $\Delta \in (0,h/2)$ \cite{lit4}.
Introduction of the shift $\Delta$ enables   singularity overcome when it comes to 
quadrature rules application. But such algorithms need this shift selection. 

\item {\it Quadrature type methods.}
One of the popular methods   (here readers may refer to work \cite{lit3})
is  Discrete Vortices Method where the integral with Cauchy kernel
$$\frac{1}{2\pi} \int\limits_{-1}^1 \frac{\gamma(x)}{x-x_0}\, dx = f(x_0), \, -1<x_0<1, $$
is approximated with lift rectangles quadrature rule and using meshes on 
$x$ and $x_0$ with shift $\Delta =h/2.$
This gives the system of linear algebraic equations (SLAE) with non zero main diagonal.

In work \cite{lit8} the ``onion peeling'' method for solution of SIE \eqref{eq1} is suggested.
Here region $r \in [0,R]$ is approximated with rings $\Delta r$ wide of 
constant values $k \in (r_j-\Delta r/2, r_j+\Delta r/2 )$ for each $r_j.$
Here meshes are assumed to be uniform ($\Delta =0$). The main idea in this method
is that integral $\int_{r_j-\Delta r/2}^{r_j+\Delta r/2} \frac{r}{\sqrt{r^2 - x_i^2}}\, dr, r_j-\Delta r/2 \geq x_i$
is computed analytically and its finite. Further midpoint quadrature is used resulting
systems of linear algebraic equations with upper triangular matrix  with respect to  $k_j = k(x_j).$
The similar method is suggested in \cite{lit37}.

\item {\it Solution approximation.}
In works \cite{lit2,lit24,lit30,lit31,lit34,lit40} et al.,
the desired solution $k(r)$ (as well as the right-hand side $q(x)$) is approximated with an
orthogonal polynomial, shifted Legendre polynomials, normalized Bernstein polynomials,
algebraic or trigonometric polynomial or polynomial spline
with coefficients determined with minimum of discrepancy between the left-hand side
and right-hand side of \eqref{eq1}. This leads to a projection method (the Galerkin method,
the collocation method, the method of splines, the quadrature method, the least squares
method, etc.) and to the solution of a SLAE wrt the corresponding polynomial coefficients.

In these algorithms, there is a {\it self-regularization}, and in the case of using the
relative shift of meshes, the shift $\Delta$ plays the role of the regularization parameter.
Namely if $\Delta$ is closer to $h/2$ then solution $k(r)$ is more stable, but it makes
reduction of resolving capability of the method. If $\Delta$ is closer to zero, then
solution is less stable but resolving capability of the method is higher.
In all these algorithms, a SLAE is with prevailing (but not infinite) matrix diagonal.

We also note a number of algorithms. Equation \eqref{eq1}, as is known, has an analytical solution
\cite{lit1,lit2,lit6,lit8,lit24,lit30,lit37}
\begin{equation}
k(r) = -\frac{1}{\pi} \int\limits_r^R \frac{q^{\prime}(x)}{\sqrt{x^2 - r^2}} \, dx, \quad 0\leq r \leq R.
\label{eq3}
\end{equation}
However, solution \eqref{eq3} contains derivative $q^{\prime}(x)$ of experimental (noisy) function
$q(x)$ and the problem of differentiation is ill-posed \cite{lit38}. Moreover, integral in \eqref{eq3}
is improper (singular). Nevertheless, a number of the following algorithms is proposed to compute
the solution according to \eqref{eq3}.

\item {\it Interpolation and quadrature method.}
In \cite{lit6},  derivative $q^{\prime}(x)$ was computed using
interpolation on three (and two) neighboring points (discrete values of $x$).
Integral $\int_r^R \frac{dx}{\sqrt{x^2 - r^2}}$ (cf. \eqref{eq3})
is computed analytically (without singularity). The similar algorithm was
suggested in \cite{lit37} using generalized left rectangles formula.

\item{\it Approximation of the right-hand side $q(x)$} is used in works \cite{lit18,lit24,lit40}.
Function $q(x)$ is suggested to be approximated by a linear combination of
smoothing polynomials (or splines) uniform for the whole interval $x \in [0,R]$.
Derivative $q^{\prime}(x)$ is computed using polynomial (or spline) differentiation.
Solution $k(r)$ in accordance with \eqref{eq3} is computed by summing the integral
in \eqref{eq3} along segments that performed analytically (see \cite[pp. 188--189]{lit40}).

\item{\it Algorithm without using derivative $q^{\prime}(x)$.}
In \cite{lit9}, formula \eqref{eq3} is converted (by means of integration by parts)
into the following expression that does not contain derivative $q^{\prime}(x)$
(cf. \cite{lit8,lit40}):
$$
k(r) = -\frac{1}{\pi} \left\{ \frac{q(R)-q(r)}{\sqrt{R^2-r^2}} + \int\limits_r^R
\frac{x\,[q(x)-q(r)]}{\sqrt{(x^2-r^2)^3}}\, dx \right\}, \quad 0\leq r\leq R.
$$
This algorithm is implemented, e.g., in paper \cite[pp. 217--220]{lit40}
using the cubic spline (see \cite[p. 273]{lit23}) for $q(x)$.

\item{\it Use of regularization.}
Abel's equation \eqref{eq1} enjoys self-regularizing property due to the
singularity, as a result the problem of its solving is moderately ill-posed \cite{lit8}.
This means that above mentioned algorithms are moderately stable.
Nevertheless, in papers \cite{lit1,lit7,lit8,lit14} et al., the Tikhonov regularization
method \cite{lit10,lit16,lit38,lit139} was used to enhance the stability of algorithms.

\end{enumerate}

In this work, we develop the following variant for numerical solving some SIE.
We make the meshes of nodes in $r$ and $x$ coincide (i.e. $\Delta=0$) and eliminate the
singularities using the generalized quadrature formula (cf. \cite{lit19,lit36,lit37}).
However, such a technique can be applied to not all SIE.
For example, it is not applicable to SIE with the Cauchy kernel, but it is applicable to
some SIE with logarithmic and other (weakly singular) kernels.
In this paper, we consider the solution of equation \eqref{eq1} wrt the desired
function $k(r)$, as well as numerical computation of $k(r)$ according to \eqref{eq3} by
the generalized quadrature method with use of Tikhonov regularization.

\section{The generalized quadrature method}

Let us describe method using generalized left rectangles formula
in application to numerical solving SIE \eqref{eq1} (the first method)
and to computation of $k(r)$ according to \eqref{eq3} (the second method).

It is to be noted here that in \cite{lit6,lit8}, the numerical method ``onion-peeling''
is suggested for computation of integrals in \eqref{eq1} \cite{lit8} and in \eqref{eq3}
\cite{lit6}. Here, the uniform coinciding node meshes in $r$ and $x$ and the middle
rectangles quadrature formula have been employed.
In \cite{lit37}, also uniform coinciding meshes have been employed combined with
more usable the left rectangles formula.

In present paper, we use nonuniform meshes and left rectangles
resulting more generic and convenient algorithm.
The solution error estimates for equation \eqref{eq1} by the generalized quadrature
method are also derived. This method is described below in two variants (the first
and second methods).

\subsection{First quadrature method}

First quadrature method employs generalized left rectangle formula.
Let us introduce nonuniform (but coinciding) meshes on $x$ and $r$ as follows
\begin{equation}
0=x_1 =r_1 <x_2 =r_2 < \ldots < x_i=r_i< \ldots <x_n =r_n =R.
\label{eq4}
\end{equation}
Here, $R=r_{\max}$ is boundary value such as $k(R+0)=0$.
On each interval $[r_j, r_{j+1})$, $j={1,2,\dots,n-1}$ we suppose approximately
\begin{equation}
k(r) = k(r_j) \equiv k_j = \text{const}.
\label{eq5}
\end{equation}


We have the following
\begin{lem}
Under condition \eqref{eq5}, one has the equality
\begin{equation}
\aligned
\int\limits_{r_j}^{r_{j+1}} \frac{r}{\sqrt{r^2-x^2}}\, k(r)\, dr =
\left(\sqrt{r_{j+1}^2 - x^2} - \sqrt{r_{j}^2 - x^2}\,\right) k_j, \\
j \in [1,n-1], \quad x\leq r_j <r_{j+1} \leq R.
\endaligned
\label{eq6}
\end{equation}
\end{lem}

\noindent {\bf Proof.} Integral (table)
$$\int\frac{r}{\sqrt{r^2-x^2}}\,dr = \sqrt{r^2-x^2} \quad \text{for} \quad x \leq r,$$
whence, taking into account (5), we obtain (6). \qquad $\Box$

\begin{definition}
We call formula \eqref{eq6} as {\it generalized quadrature formula of left rectangles}
(cf. \cite{lit19}) for the specific  singularity $r / \sqrt{r^2 - x^2}$, and
multipliers $\sqrt{r_{j+1}^2 - x^2} - \sqrt{r_{j}^2 - x^2} $ are
the {\it quadrature coefficients } of this singularity.
\end{definition}

The specifics of the formula \eqref{eq6} is that the singular integral
$\int_{r_j}^{r_{j+1}}\frac{r}{\sqrt{r^2-x^2}}\,dr$
is calculated analytically accurate and without peculiarity.
If it is calculated numerically by the usual left rectangles quadrature formula,
then at $x = r_j$ there will be a division by zero.

Let us now formulate the main result as following

\begin{theorem}
Numerical solution of equation \eqref{eq1} according to the first
quadrature method is defined as following  recursion
\begin{equation}
\left\{
\begin{array}{l}
k_{n-1} = \frac{\displaystyle q_{n-1}/2}{\displaystyle p_{n-1,n-1}}, \\[10pt]
k_i= \frac{\displaystyle q_i/2 - \sum\nolimits_{j=i+1}^{n-1} p_{ij}k_j}{\displaystyle p_{ij}},
\,\,\, i=n-2,n-3, \dots, 1, \\[10pt]
k_n =k_{n-2} + \left(\frac{\displaystyle r_n - r_{n-2}}{\displaystyle r_{n-1} - r_{n-2}}
\right)(k_{n-1}-k_{n-2}),
\end{array}
\right.
\label{eq7}
\end{equation}
where $k_i \equiv k(r_i)$, $q_i \equiv q(x_i)$,
\begin{equation}
p_{ij}=\sqrt{r_{j+1}^2 - x_i^2} - \sqrt{r_{j}^2 - x_i^2},
\label{eq8}
\end{equation}
\end{theorem}

\noindent {\bf Proof.}
Integral in \eqref{eq1} is sum of integrals \eqref{eq6}, i.e.
\begin{equation}
\aligned
\int\limits_{x_i}^{R} \frac{r}{\sqrt{r^2-x_i^2}} \, k(r) \, dr = \sum\limits_{j=i}^{n-1}
\left(\sqrt{r_{j+1}^2 - x_i^2} - \sqrt{r_{j}^2 - x_i^2}\,\right) k_j = q_i/2, \\
i = {1,2,\dots,n-1}.
\endaligned
\label{eq9}
\end{equation}
This is the SLAE wrt $\{k_j\}_{j=1}^n$.
SLAE \eqref{eq9} is upper triangular and its solution can be recursively
constructed. From \eqref{eq9} for $i=n-1,n-2, \dots, 1$
we eventually obtain $k_{n-1},\, k_{n-2}, \dots , k_1$ according to \eqref{eq7}.
As to the value of $k_n \equiv k(R)$, it can not be found by this scheme, but can be
additionally determined as $k_n = 0$ from physical concepts or $k_n = k_{n-1}$ or
can be derived using linear extrapolation \cite{lit35}, as was done in \eqref{eq7}.
\qquad $\Box$

In \cite{lit8}, formulae of type \eqref{eq6} are also given, but for the case of uniform
(and coinciding) meshes in $r$ and $x$ and using the middle rectangles quadrature formula.
Furthermore, important formulae of type \eqref{eq7} is not given.

Formulae \eqref{eq4}--\eqref{eq9} are more common and more convenient than in \cite{lit8}
and formulae \eqref{eq7} give a solution in the explicit form. The method according to
\eqref{eq4}--\eqref{eq9} for solving the equation \eqref{eq1} is called in \cite{lit37}
the {\it generalized quadrature method} for solving SIE \eqref{eq1}. This paper presents
a more general formulae \eqref{eq4} and \eqref{eq7} than in \cite{lit37}.
In Sec. 3, estimates of the errors for this method are given.


\subsection{Second  quadrature method}

The second method is generalized quadrature method for computation of singular integral \eqref{eq3}
giving the solution $k(r)$. Let us assume $k^{\prime}(x)$ to be computed with some
stable method. As in the first method we introduce node meshes \eqref{eq4}.
On each $[x_i, x_{i+1})$, $i={1,2, \dots, n-1}$ we assume
\begin{equation}
q^{\prime}(x)=q^{\prime}(x_i)\equiv q_i^{\prime} = \text{const}.
\label{eq10}
\end{equation}
Let us fomulate the following
%
\begin{lem}
Under condition \eqref{eq10} the following equality is true
\begin{equation}
\aligned
\int\limits_{x_i}^{x_{i+1}} \frac{q^{\prime}(x)}{\sqrt{x^2-r^2}} \, dx =
\ln\frac{x_{i+1}+\sqrt{x_{i+1}^2-r^2}}{x_{i}+\sqrt{x_{i}^2-r^2}}\,q_i^{\prime},\\
i=1,2,\dots,n-1, \quad r\leq x_i<x_{i+1} \leq R.
\endaligned
\label{eq11}
\end{equation}
\end{lem}

\noindent {\bf Proof.}
Integral
$\int\frac{dx}{\sqrt{x^2-r^2}} = \ln\left(x+\sqrt{x^2-r^2}\,\right)$ for $r\leq x$
is the table integral.
Taking into account the condition \eqref{eq10} we obtain \eqref{eq11}. \qquad $\Box$

\begin{definition}
Formula \eqref{eq11} is {\it generalized left rectangle quadrature rule }
for the singularity $1/\sqrt{r^2-x^2}$, and multipliers
$\ln\frac{x_{i+1}+\sqrt{x_{i+1}^2-r^2}}{x_{i}+\sqrt{x_{i}^2-r^2}}$
are {\it quadrature coefficients} of this singularity.
\end{definition}

Now we can formulate
\begin{theorem}
Numerical solution of SIE \eqref{eq1} by formula \eqref{eq3} according to
the second generalized quadrature method
is result of the following recurrence formulae
\begin{equation}
\left\{
\begin{array}{l}
k_j = -\frac{\displaystyle 1}{\displaystyle \pi}\displaystyle\sum\limits_{i=j}^{n-1}
g_{ij}\,q_i^{\prime}, \quad j={2,3, \dots, n-1}, \\[15pt]
k_1=\frac{\displaystyle q_1/2 - \sum\nolimits_{j=2}^{n-1}(r_{j+1}-r_j)\,k_j}{\displaystyle r_2},
\end{array}
\right.
\label{eq12}
\end{equation}

%
where
\begin{equation}
g_{ij} = \ln \frac{x_{i+1}+\sqrt{x_{i+1}^2-r_j^2}}{x_{i}+\sqrt{x_{i}^2-r_j^2}}.
\label{eq13}
\end{equation}
\end{theorem}

\noindent {\bf Proof.}
Integral in \eqref{eq3} is sum of integrals \eqref{eq11} over the separate intervals
$[x_i,x_{i+1})$, i.e.
\begin{equation}
\int\limits_{r_j}^{R} \frac{q^{\prime}(x)}{\sqrt{x^2-r_j^2}} \, dx =
\sum_{i=j}^{n-2}\ln\frac{x_{i+1}+\sqrt{x_{i+1}^2-r_j^2}}{x_{i}+\sqrt{x_{i}^2-r_j^2}}\,
q_i^{\prime}, \quad j={1,2, \dots, n-1}.
\label{eq14}
\end{equation}
As a result, solution \eqref{eq3} in the discrete form is $\{ k_j \}_{j=2}^{n-1}$
according to \eqref{eq12}, \eqref{eq13}.
For $j=1$ this (second) method due to \eqref{eq14} gives uncertainty $\infty \cdot 0$
for $i=j=1$ since $x_1=r_1=q_1^{\prime}=0.$ In this case let's use the first method to
determine $k_1$. Using \eqref{eq7} and \eqref{eq8} for $i=1$ we find $k_1$, ref. \eqref{eq12}.
As to $k_n$ it can be calculated using linear extrapolation (see \eqref{eq7}). \qquad $\Box$

The advantage of the above two methods is that they do not require the relative shift of
meshes and their integrals with singularities $r/\sqrt{r^2-x^2}$ and $1/\sqrt{r^2-x^2}$ are
calculated analytically and without divergences.

However, these methods are not suitable for all singularities, e.g., for numerical
computation of hypersingular integral with the Cauchy kernel
$\int_{-1}^1 \frac{x(\tau)}{\tau -t} \, d\tau$ \cite{lit3,lit4}.


%
%

\section{Error estimates}

Let us derive errors estimates for solution of SIE \eqref{eq1} using fist quadrature method
(cf. \cite{lit37,lit39,lit20}).

\subsection{Quadrature error on small interval}

Let us estimate quadrature error for computing integral \eqref{eq6} on separate small interval
$[r_j, r_{j+1})$ due to approximation \eqref{eq5} (while without the measurement error for $q(x)$).

\begin{lem}
Integral
\begin{equation}
\aligned
\int\limits_{r_j}^{r_{j+1}} \frac{r}{\sqrt{r^2-x_i^2}}\,k(r)\,dr,
\quad x_i\leq r_j < r_{j+1}\leq R,\\
i={1,2,\dots,n-1}, \quad j={i,\dots, n-1},
\endaligned
\label{eq15}
\end{equation}
when using the generalized left rectangles formula \eqref{eq6} and taking account of
quadrature error caused by the approximation \eqref{eq5} is equal to (refinement of
formula \eqref{eq6})
\begin{equation}
\int\limits_{r_j}^{r_{j+1}} \frac{r}{\sqrt{r^2-x_i^2}}\,k(r)\,dr=
p_{ij}k_j+\Delta\varepsilon_{ij},
\label{eq16}
\end{equation}
where $p_{ij}$ are quadrature coefficients \eqref{eq8} and $\Delta\varepsilon_{ij}$ is
quadrature error of computation of integral \eqref{eq15} approximately equal
\begin{equation}
\begin{array}{rcl}
\Delta\varepsilon_{ij} &=& \frac{\displaystyle k^{\prime}(\xi_j)}{\displaystyle 2} \Biggl[
(r_{j+1} - 2r_j)\sqrt{r_{j+1}^2-x_i^2}+r_j\sqrt{r_{j}^2-x_i^2} \\
&+& x_i^2\ln\frac{\textstyle r_{j+1}+\sqrt{r_{j+1}^2-x_i^2}}
{\textstyle r_{j}+\sqrt{r_{j}^2-x_i^2}}\, \Biggr],
\quad \xi_j \in [r_j, r_{j+1}).
\end{array}
\label{eq17}
\end{equation}
\end{lem}

\noindent {\bf Proof.}
Using the first method we assume $\widetilde{k}(r) = k_j$, $r\in[r_j,r_{j+1})$ (see \eqref{eq5}),
i.e. we represent function $k(r)$ by the interpolation Lagrange zero degree polynomial \cite{lit17}.
Error of such interpolation is $\Delta k_j(r) \equiv k(r) - k_j = k^{\prime}(\xi)\,(r-r_j)$,
where $\xi=\xi_j(r)\in[r_j,r_{j+1})$. Then
\begin{equation}
k(r) = k_j +k^{\prime}(\xi)\,(r-r_j).
\label{eq18}
\end{equation}
Let us now substitute \eqref{eq18} into \eqref{eq15}, we get \eqref{eq16}, where
\begin{equation}
\Delta\varepsilon_{ij} = k^{\prime}(\xi_j) \int\limits_{r_j}^{r_{j+1}}
\frac{r\,(r-r_j)}{\sqrt{r^2 - x_i^2}} \, dr.
\label{eq19}
\end{equation}
Integral in \eqref{eq19} can be analytically computed giving us an estimate \eqref{eq17}.
\qquad $\Box$

It is to be noted that derivative $k^{\prime}(\xi_j)$ in \eqref{eq17} can be approximated
with
\begin{equation}
k^{\prime}(\xi_j) = \frac{k_{j+1} - k_j}{r_{j+1}-r_j}, \quad j={1,2,\dots,n-1},
\label{eq20}
\end{equation}
or by other way \cite{lit37}.

\subsection{Quadrature error}

Let us estimate quadrature error of the solution of equation \eqref{eq1} due to the approach
\eqref{eq5} (without error $q(x)$). We formulate this as a theorem.

\begin{theorem}
Errors of numerical solution of equation \eqref{eq1} by the first
generalized quadrature method according to \eqref{eq7} are computed
with the following recurrence
\begin{equation}
\left\{ \begin{array}{l}
\Delta k_n = \Delta k_{n-1} = \frac{\displaystyle\Delta \varepsilon_{n-1,n-1}}
{\displaystyle p_{n-1,n-1}}, \\[10pt]
\Delta k_i = \frac{\displaystyle\varepsilon_i - \sum\nolimits_{j=i+1}^{n-1} p_{ij} \,
\Delta k_j}{\displaystyle p_{ii}}, \quad i=n-2, n-3, \dots, 1,
\end{array}
\right.
\label{eq21}
\end{equation}
where
\begin{equation}
\varepsilon_i = \sum_{j=i}^{n-1} \Delta \varepsilon_{ij}.
\label{eq22}
\end{equation}
Here, $p_{ij}$, $\Delta\varepsilon_{ij}$ and $k^{\prime}(\xi_j)$
are computed based on \eqref{eq8}, \eqref{eq17} and \eqref{eq20} respectively.
\end{theorem}

\noindent {\bf Proof.}
Let us write integral in \eqref{eq1} as sum of integrals wrt intervals
$$
\int\limits_{x_i}^R \frac{r}{\sqrt{r^2 -x_i^2}} \, k(r) \, dr = \sum_{j=i}^{n-1}
\int\limits_{r_j}^{r_{j+1}} \frac{r}{\sqrt{r^2 - x_i^2}} \, (k_j +\Delta k_j(r)) \, dr
$$
$$
=\sum_{j=i}^{n-1} \int\limits_{r_j}^{r_{j+1}} \frac{r}{\sqrt{r^2-x_i^2}} \, (k_j +
k^{\prime}(\xi)\,(r-r_j))\,dr, \quad i={1,2,\dots,n-1}.
$$
Then
\begin{equation}
\sum_{j=i}^{n-1} \int\limits_{r_j}^{r_{j+1}} \frac{r}{\sqrt{r^2 - x_i^2}}\,\Delta k_j(r)\,dr
= \sum_{j=i}^{n-1} \int\limits_{r_j}^{r_{j+1}} \frac{r\,(r-r_j)}{\sqrt{r^2 - x_i^2}}\,k^{\prime}(\xi)\, dr,
\,\, i={1,2,\dots,n-1}.
\label{eq23}
\end{equation}

In order to compute the integral in the left-hand side of \eqref{eq23}, we employ the formula of
left rectangles, i.e. we assume $\Delta k_j(r) = \Delta k(r_j) \equiv \Delta k_j = \text{const}$,
$r\in [r_j,r_{j+1})$ and compute the integral using the generalized formula in similar way with \eqref{eq6}.
Integral in the right-hand side of \eqref{eq23} is equal to $\Delta\varepsilon_{ij}$
due to \eqref{eq19}. Then
\begin{equation}
\sum_{j=i}^{n-1} p_{ij} \, \Delta k_j = \varepsilon_i, \quad i={1,2,\dots,n-1},
\label{eq24}
\end{equation}
where $\varepsilon_i$ denote sum \eqref{eq22}. Here,
\eqref{eq24} is a SLAE wrt $\{ \Delta k_j \}_{j=1}^{n-1}$.
It is also assumed that $\{k_i \}_{i=1}^{n}$ are computed using \eqref{eq7} in advance.
Then values $\{ \varepsilon_i  \}_{i=1}^{n-1}$ are computed using
\eqref{eq22}, \eqref{eq17} and \eqref{eq20}. We solve SLAE \eqref{eq24}
with upper triangular matrix and obtain solution \eqref{eq21} in the recurrent form,
adding the condition $\Delta k_n = \Delta k_{n-1}$. \qquad $\Box$

\noindent {\bf Remark 2.} It is to be noted here that formulae \eqref{eq21} give errors of the solution $\Delta k_i$
with their signs (cf. \cite{lit37,lit39}) in contrast with other works where absolute values $|\Delta k_i|$
or upper bounds $|\Delta k_i|\leq\dots$ or upper bounds by the norm $\|\Delta k_i\|\leq\dots$, etc. are given.
This enable us to obtain the refined solution
\begin{equation}
\widehat{k}_i = k_i +\Delta k_i, \quad i={1,2,\dots,n},
\label{eq25}
\end{equation}
using $\{k_i \}_{i=1}^n$ from \eqref{eq7} and errors $\{\Delta k_i \}_{i=1}^n$ from \eqref{eq21}.

\noindent {\bf Remark 3.}
Errors of numerical solution given in \eqref{eq21} are obtained with regard to only the quadrature
errors and the error of the right-hand side $q(x)$ of equation \eqref{eq1} is set equal to zero.
Let us take into account the measurement errors $\{\delta_i\}_{i=1}^{n}$ of the source function $q(x)$.
In \cite{lit20,lit39}, the error estimates for numerical solution of the Volterra integral
equations of the first and second kind are derived taking into account both quadrature and
source function errors. In similar way we can generalize recurrence formulae \eqref{eq21}
to the case of errors $\{\delta_i\}_{i=1}^{n}$ as follows
\begin{equation}
\left\{ \begin{array}{l}
|\Delta k_n| = |\Delta k_{n-1}| = \frac{\displaystyle|\Delta\varepsilon_{n-1,n-1}| +
\delta_{n-1}}{\displaystyle p_{n-1,n-1}}, \\[10pt]
|\Delta k_i| = \frac{\displaystyle|\varepsilon_i|+\delta_i+\sum\nolimits_{j=i+1}^{n-1}
p_{ij}\,|\Delta k_j|}{\displaystyle p_{ii}}, \quad i=n-2,n-3,\dots,1.
\end{array}
\right.
\label{eq26}
\end{equation}
But estimates \eqref{eq26} give overstated estimates of $\{|\Delta k_i|\}_{i=1}^{n}$ due to
using the operation $|\cdot|$ (absolute value).

\section{Numerical illustration}

\subsection{Software implementation}

Proposed two generalized quadrature methods have been implemented in MatLab\,7.10 (R2010a).
Following the first method we search $\{k_i\}_{i=1}^n$ using \eqref{eq7}, \eqref{eq8}.
Errors $\{\Delta k_i\}_{i=1}^n$ are calculated using \eqref{eq21}, as well as \eqref{eq8}, \eqref{eq17},
\eqref{eq20} and \eqref{eq22}. Refined solution $\{\widehat{k}_i\}_{i=1}^n$ is calculated with \eqref{eq25}.
Tikhonov regularization  \cite{lit10,lit16,lit38,lit39}
\begin{equation}
k_{\alpha} = (\alpha E + A^T A)^{-1} A^T f
\label{eq27}
\end{equation}
is employed, where the discrepancy principle \cite{lit25} is used for choosing
the regularization parameter $\alpha>0$:
\begin{equation}
\|Ak_{\alpha}-f\| = \delta.
\label{eq28}
\end{equation}
Here, $f=q/2$, $E$ is identity matrix, $A$ is matrix of the SLAE \eqref{eq9} represented as
\begin{equation}
Ak=f,
\label{eq29}
\end{equation}
where
\begin{equation}
A_{ij}=\left\{
\begin{array}{ll}
p_{ij}, &  j\geq i, \\
0, &  \text{otherwise},
\end{array}
\right.
\quad i,j=1,2,\dots,n-1.
\label{eq30}
\end{equation}

Following the second method, solution $k(r)$ has been  computed using singular integral in
\eqref{eq3} based on generalized left rectangles formula according to \eqref{eq12} and \eqref{eq13}.

The first method was developed and implemented in software in more detail than the second method.

\subsection{Infrared tomography example}
\begin{figure}[htp] 
\centering\includegraphics[width=9cm]{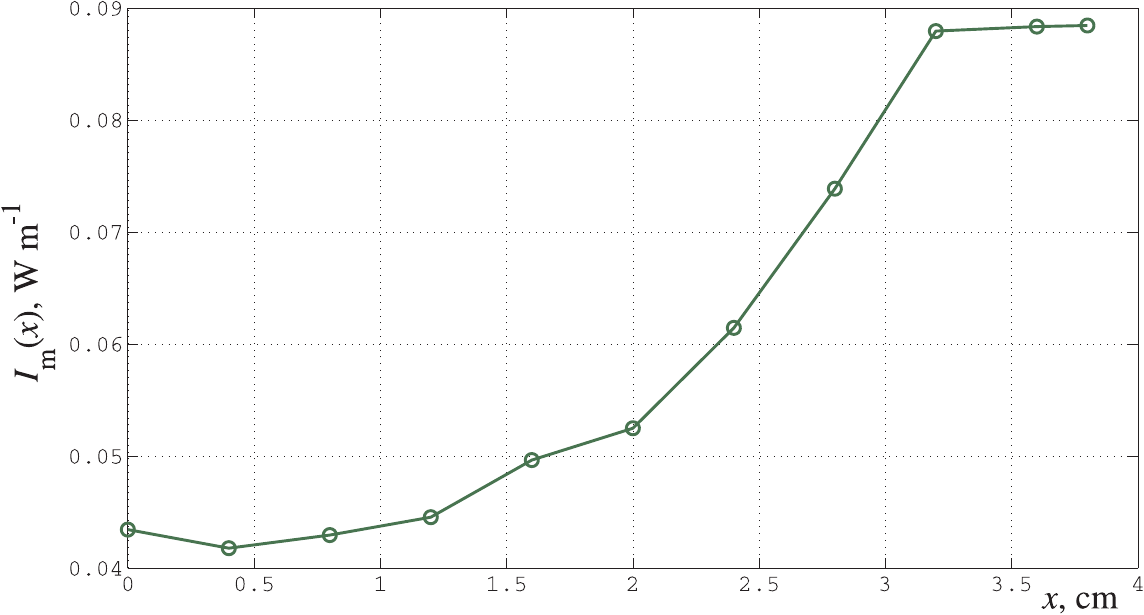} 
\caption{Measured noisy intensity $I_{\text{m}}(x)$ (difference of intensity in active and
passive regimes). The mesh is nonuniform, number of nodes $n=11$.}
\label{fig1}
\end{figure}

The first method has been applied to axially symmetric flame diagnostics using infrared tomography
\cite{lit1,lit6,lit7,lit8,lit11,lit12,lit15,lit29,lit36}.
Fig. \ref{fig1} shows measured output intensity $I_{\text{m}}(x)$ of rays (m from measurement)
which go through gas, undergo absorption and emission and are accepted by detectors.
Measurements have been performed in Technical University of Denmark, Department of Chemical and
Biochemical Engineering (before 1 January 2012 Ris$\o$ DTU) within a joint project \cite{lit12,lit11}.

\begin{figure}[htp] 
\centering{\includegraphics[width=9cm]{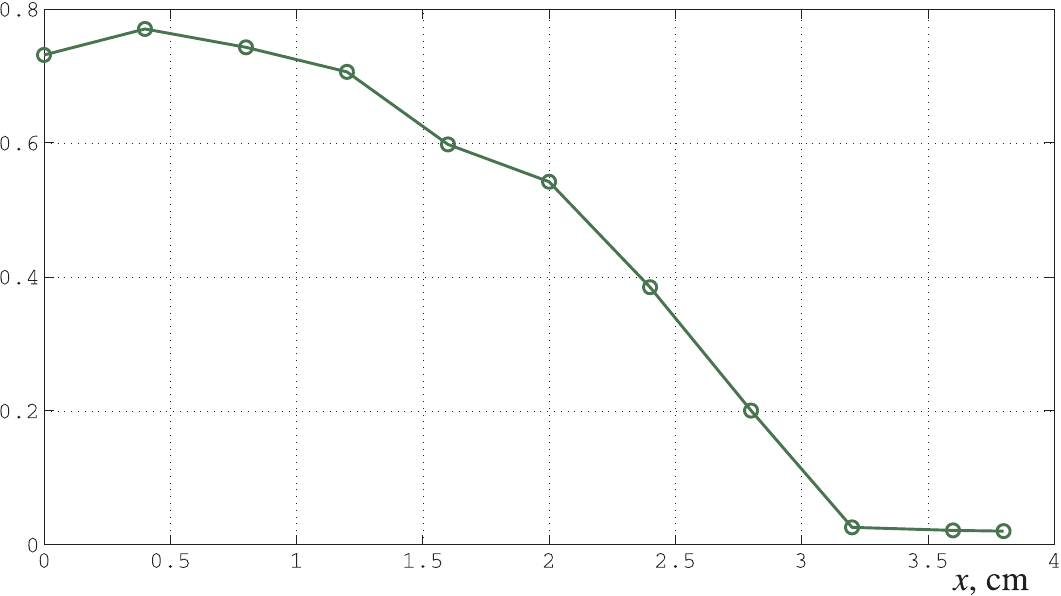}} 
\caption{Dimensionless right-hand side $q_{\text{m}}(x)$ of SIE \eqref{eq1}, $n=11$.}
\label{fig2}
\end{figure}

Intensity $I_{\text{m}}(x)$ is recalculated into $q_{\text{m}}(x) = -\ln[I_{\text{m}}(x)/B(T_0)]$
(the right-hand side of equation \eqref{eq1}), where $B(T_0)$ is the Planck function of rays
source with its temperature $T_0=894.4^{\circ}$C. Fig. \ref{fig2} shows function $q_{\text{m}}(x)$.

\begin{figure}[htp] 
\centering{\includegraphics[width=9cm]{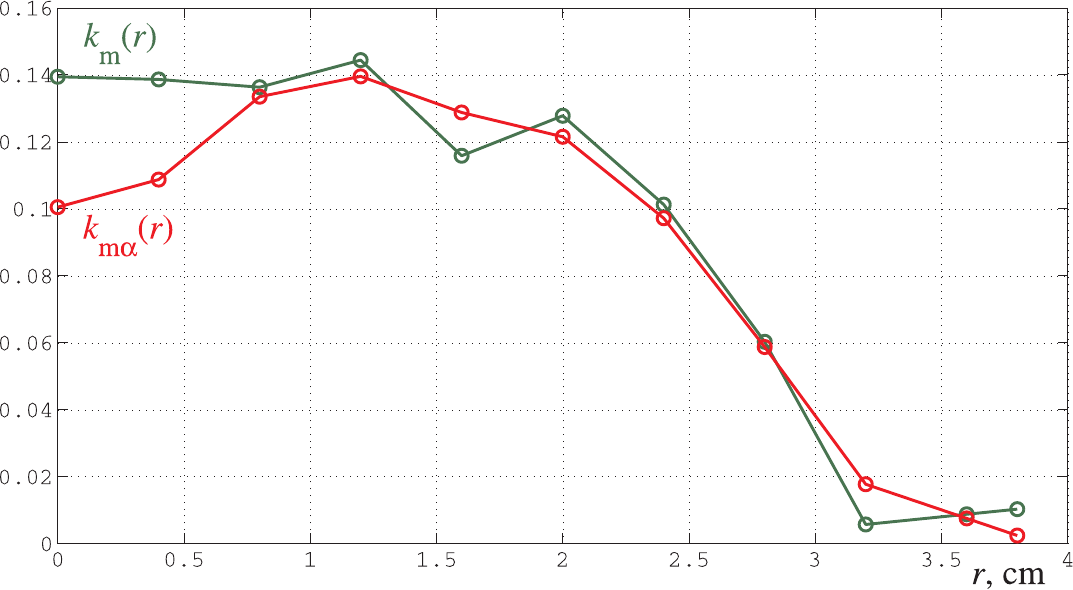}} 
\caption{Absorption coefficient $k_{\text{m}}(r)$ computed by the first generalized quadrature
method and $k_{\text{m}\alpha}(r)$ computed by Tikhonov regularization, $\text{cm}^{-1}$. }
\label{fig3}
\end{figure}

Fig. \ref{fig3} shows results of solution of SIE \eqref{eq1} wrt absorption coefficient
$k_{\text{m}}(r)$ by the first method of generalized quadratures according to \eqref{eq7}
and \eqref{eq8}.


As we see, the solution $k_{\text{m}}(r)$ suffer from significant artificial perturbations.
This is due to too great step of node mesh (in other words, the smallness of $n$), as well as
measurement errors in function $I_{\text{m}}(x)$.
Here we also demonstrate  behavior of solution $k_{\text{m}\alpha}(r)$ derived with Tikhonov
regularization using \eqref{eq27}, \eqref{eq29}, \eqref{eq30} and \eqref{eq9}.
Regularization parameter $\alpha$ is chosen  using discrepancy principle \eqref{eq28}
where $\delta = 0.037$, as a result $\alpha = 10^{-0.09} = 0.813$.
Fig. \ref{fig3} demonstrates that solution has been smoothed by regularization method.

To reduce the grid step in $x$ as well as to moderately smooth the fluctuations in the
function $I_{\text{m}}(x)$, a spline approximation was used \cite{lit23,lit35,lit40}.
Fig. 4 shows an approximation of the function $I_{\text{m}}(x)$ by cubic smoothing spline
using the m-function csaps.m.

\begin{figure}[htp] 
\centering{\includegraphics[width=9cm]{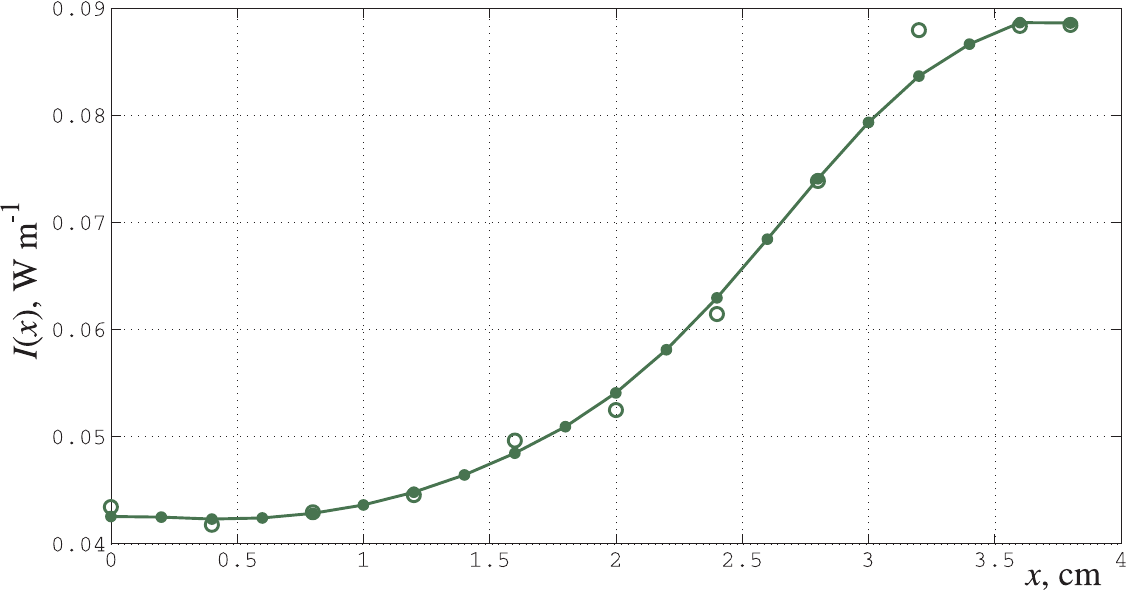}} 
\caption{Measured $I_{\text{m}}(x)$ values are marked with $\circ$ (n=11), values
$I_{\text{m}}(x)$ for $n=20$ marked with $\bullet$ and spline interpolated values
are marked with solid line. }
\label{fig4}
\end{figure}

\begin{figure}[htp] 
\centering{\includegraphics[width=9cm]{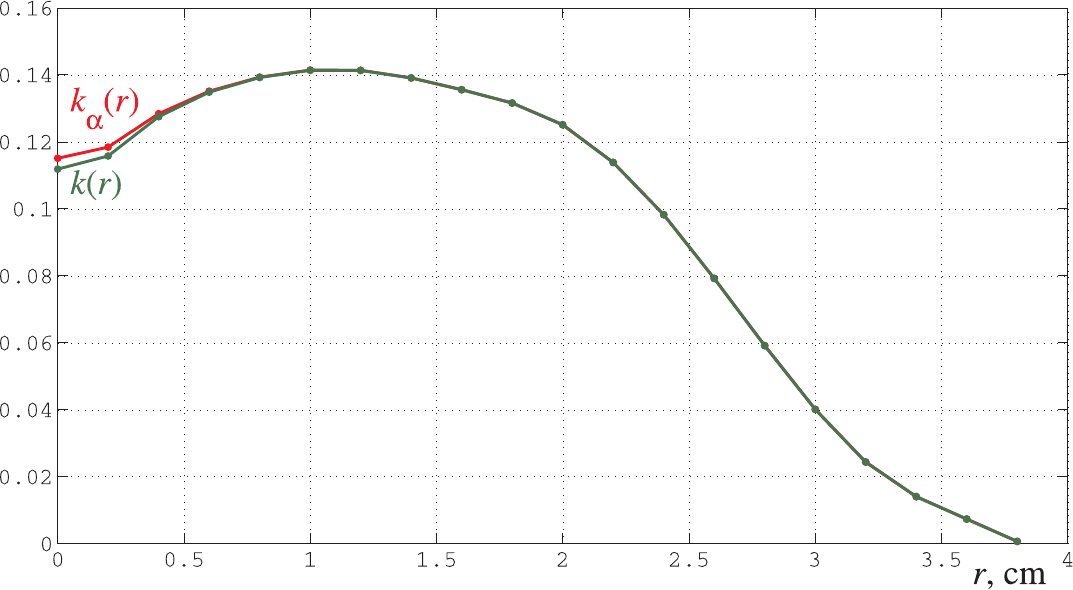}} 
\caption{Absorption coefficient $k(r)$ (without regularization) and $k_{\alpha}(r)$
(with regularization) after spline smoothing of $I_{\text{m}}(x)$, $n=20$,
$\text{cm}^{-1}$. }
\label{fig5}
\end{figure}

The smoothed values of $I(x)$  were generated with splines (Fig. \ref{fig4})
and then  SIE \eqref{eq1} was resolved with generalized quadrature \eqref{eq7}.
Fig. \ref{fig5} shows obtained solution $k(r)$.
We also applied solution using Tikhonov regularization \eqref{eq27} for
$\alpha = 10^{-2}$. Fig. \eqref{fig5} shows regularized solution $k_\alpha(r)$.

Fig. \ref{fig4} and \ref{fig5} demonstrate that application of spline based smoothing
enable mesh step reduction for $x$ (causing increase of $n$). This allows (moderate)
smoothing $k(r)$ and $k_{\alpha}(r)$. Moreover, in case of noisy $I(x)$
and big step of the mesh, regularization slightly improves solution as shown in Fig. \ref{fig3}.
In case of $< 1$\% errors and small step (Fig. \ref{fig4}) solutions $k(r)$ (without
regularization) and $k_{\alpha}(r)$ (with regularization) are obtained practically
the same (Fig.\ref{fig5}).
It confirms that the problem of solving the singular integral equations is moderately
ill-posed and has the property of self-regularization.

  \section{Conclusion}

In this paper we outlined two new methods of numerical solution of singular integral equation (SIE)
of Abel type. The methods are based on the use of generalized quadrature formula of left rectangles.
Specificity of methods is that singular integrals are computed analytically and without peculiarities.
We derive recurrence formulae for solution, generally speaking, on a nonuniform node mesh.
Estimates of quadrature errors of solution with regard to their sign in the absence and in the presence
of errors in the right-hand side are found.
In order to enhance the stability of the solution we used Tikhonov regularization.
However, SIE enjoy self-regularization, therefore it is advisable to apply the Tikhonov regularization
method only if there is a significant error ($\sim10\%$) in the right-hand side and rough mesh step
(when the number of nodes is small: $n\sim10$).
The method has been applied for solution of infrared tomography problem.

\section*{Asknowledgements}

This work is supported with RFBR (Projects No. 09-08-00034,  No. 13-08-00442) and DTU, Denmark
(Project No. 010246).
The authors thank V. Evseev and A. Fateev for experimental data (Fig. \ref{fig1}) and for useful
discussions.





\bibliographystyle{elsarticle-num}


 

\end{document}